\documentclass[11pt]{amsproc}

\usepackage{amssymb}

\theoremstyle{plain}
\newtheorem{thm}{Theorem}[section]
\newtheorem{lem}[thm]{Lemma}

\theoremstyle{remark}
\newtheorem{rem}[thm]{Remark}

\begin{document}
\title[]{On a homotopy equivalence between the
$2$-local geometry and the Bouc complex for the sporadic group
$\mathbf{McL}$}
\maketitle

\begin{center}
{\Large John Maginnis\footnote{Email address: maginnis@math.ksu.edu.} and Silvia Onofrei\footnote{Email address: onofrei@math.ksu.edu.}}\\
{\it Department of Mathematics, Kansas State University,\\ 138 Cardwell
Hall, Manhattan, KS 66506}

\medskip
{\it Archiv der Mathematik (to appear)}\\
\end{center}
\begin{abstract}We study the homotopy relation between the
standard $2$-local geometry $\Delta$ and the Bouc complex for the sporadic
finite simple group $McL$.

\smallskip
\noindent{\it Keywords:} Subgroup complexes, $2$-radical subgroups, Bouc complex, 2-local geometry, homotopy equivalence.\\
\noindent{\it 1991 MSC}: 20G05, 20C20, 51D20.\\
\end{abstract}

\section{Introduction}

The homotopy equivalence between the pure $2$-local geometry $\Delta$ for the
sporadic simple group McL and the Quillen collection was proved in \cite{sy}. The
Quillen collection is homotopy equivalent to the Bouc collection of all
non-trivial $p$-radical subgroups. It follows that for $G=\rm{McL}$ the $2$-local
geometry $\Delta$ and the Bouc collection of non-trivial $2$-radical subgroups
are also homotopy equivalent.

\smallskip
It is the purpose of this paper to give a direct relationship bewteen the $2$-local geometry $\Delta$ and the Bouc complex $\mathcal{B}_2(\rm{McL})$. We prove that a certain subcomplex $\mathcal{B}_2^I$ of the Bouc complex, which is isomorphic to a subdivision of the geometry $\Delta$, is homotopy equivalent to the full Bouc complex $\mathcal{B}_2$ (in fact, this is a retract under an equivariant deformation retraction). Our approach is similar to the one used in our previous work \cite{mgo1}, which describes such a relationship for the sporadic group Co$_3$, although for Co$_3$ the standard $2$-local geometry is homotopy equivalent to the ``distinguished" Bouc collection; see \cite{mgo2}.

\smallskip
In Section $2$ details on the $2$-local geometry $\Delta$ are given. In Section $3$ the Bouc collection $\mathcal{B}_2(\rm{McL})$ is described using geometric interpretations for each of its subgroups. In Section $4$, the main result, Theorem $4.1$, is proved.

\medskip
{\bf Acknowledgements}
\\We would like to thank Stephen Smith for suggesting this problem to us.\\

\section{The $2$-local geometry $\Delta$ of $McL$}

Let $G$ denote the sporadic simple group of McLaughlin. This is a group with one class of involutions, of $2$-rank four and whose Sylow $2$-subgroup is of order $2^7$. The maximal subgroups of $G$ have been determined in \cite{fin73}. Some details on the subgroup structure as well the computation of the mod-$2$ cohomology of $G$ can be found in \cite{am97}.

\medskip
A sporadic geometry $\tilde \Delta$ for McL was initially given in \cite{rst}. A very nice description of this geometry can be found in \cite{str}. This is a rank $4$ geometry with the following diagram:

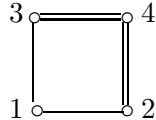
\begin{picture}(1000,50)(0,-10)
\put(145,26){$\circ$}
\put(138,27){$3$}
\put(150,28){\line(1,0){30}}
\put(150,30){\line(1,0){30}}
\put(147,-3){\line(0,1){30}}
\put(146,-9){$\circ$}
\put(138,-9){$1$}
\put(180,26){$\circ$}
\put(188,27){$4$}
\put(181,-4){\line(0,1){31}}
\put(183,-4){\line(0,1){31}}
\put(180,-9){$\circ$}
\put(188,-9){$2$}
\put(151,-7){\line(1,0){30}}
\end{picture}

\medskip
In \cite{rsy} it was shown that, for any group $P$ of order $2$ in $G$, the fixed point set $\widetilde{\Delta}^P$ is contractible. It follows that the complex $\widetilde{\Delta}$ is ample; that is, it provides a homology decomposition for the classifying space of McL. In \cite{sy} a subcomplex $\Delta$ of $\widetilde{\Delta}$ was introduced. This is a truncation of $\widetilde{\Delta}$ obtained by removing all objects of type $1$ and all the flags containing them. The new complex $\Delta$ is also ample but it has the advantage of being pure $2$-local, in the sense of \cite{sy}, which means that all the simplex stabilizers are $2$-local subgroups of $G$. The stabilizers of the objects are maximal $2$-local subgroups of $G$. Benson and Smith \cite[Section 6.11]{bs} used this complex to obtain a homology decomposition of the classifying space of $G$.

\medskip
In what follows the $2$-local geometry $\Delta$ will be used. It has diagram:

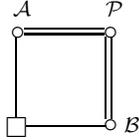
\begin{picture}(1000,50)(0,-10)
\put(145,26){$\circ$}
\put(146,35){\scriptsize{$\mathcal{A}$}}
\put(150,28){\line(1,0){30}}
\put(150,30){\line(1,0){30}}
\put(147,-3){\line(0,1){30}}
\put(132,-11){$\qed$}
\put(180,26){$\circ$}
\put(181,35){\scriptsize{$\mathcal{P}$}}
\put(181,-4){\line(0,1){31}}
\put(183,-4){\line(0,1){31}}
\put(180,-9){$\circ$}
\put(188,-9){\scriptsize{$\mathcal{B}$}}
\put(151,-7){\line(1,0){30}}
\end{picture}

The objects are denoted points $\mathcal{P}$, $\mathcal{A}$-spaces and $\mathcal{B}$-spaces. These objects correspond to the elementary abelian $2$-subgroups of $G$ as follows: the points correspond to subgroups $A_1 \simeq 2$, the $\mathcal{A}$-spaces to $A_4 \simeq 2^4$ and the $\mathcal{B}$-spaces to $A_{4'} \simeq 2^4$. Note that $A_4$ and $A_{4'}$ correspond to the two conjugacy classes of maximal elementary abelian $2$-subgroups. Incidence between elements of $\mathcal{P}$ and elements of $\mathcal{A}$ ($\mathcal{B}$ respectively) is given by containment. Two spaces $A \in \mathcal{A}$ and $B \in \mathcal{B}$ are said to be incident if the corresponding subgroups have a $2^2$ in common, that is $A_4 \cap A_{4'} = A_2$. We will refer to the conjugacy class of subgroups of type $2^2$ as lines $\mathcal{L}$. Given a subgroup $A_2 \simeq 2^2$ there exists a unique pair of incident spaces $A$ and $B$, with $A \cap B =L$; see \cite[Section 6.11]{bs}.

\medskip
The geometry $\Delta$ can be regarded as a simplicial complex
of dimension two with three types of vertices. The group $G$
acts (faithfully) flag-transitively on the geometry; see \cite{bs} for a proof. The stabilizers of
the three types of objects are:

\medskip
\hspace*{1cm}\begin{tabular}{lll}
$G_p =N_G(A_1) \simeq 2^.A_8$ && for a point $p \in \mathcal{P}$;\\
$G_A =N_G(A_4) \simeq 2^4:A_7$ && for $A \in \mathcal{A}$;\\
$G_B =N_G(A_{4'})\simeq 2^4:A_7$ && for $B \in \mathcal{B}$.\\
\end{tabular}

\medskip
The flag stabilizers can be easily determined and are as follows:

\medskip
\hspace*{1cm}\begin{tabular}{llll}
$G_{pA} \simeq 2^4:L_3(2)$ &&\phantom{aaaaaaaa}&$G_{AB} \simeq 2^4:(3 \times A_4):2$\\
$G_{pB} \simeq 2^4:L_3(2)$ &&&$G_{pAB} \simeq 2^4:S_4$.\\
\end{tabular}

\medskip
Some of the properties of  $\Delta$ which can be read from the
diagram:
\begin{enumerate}
\item[$\mathfrak{1)}.$] Given $p \in \mathcal{P}$, there are $15$ $\mathcal{A}$-spaces and $15$ $\mathcal{B}$-spaces incident with $p$; they correspond to the ``points" and the ``planes" of a projective space of rank $3$ over the field with two elements $PG(3,2)$.
\item[$\mathfrak{2)}.$]  Let $A \in \mathcal{A}$. Then there are $15$ points in $A$ and $35$ $\mathcal{B}$-spaces incident with $A$. The geometry $\text{Res}_{\Delta}A = ( \mathcal{P}_A, \mathcal{B}_A)$, of all points and all $\mathcal{B}$-spaces incident with $A$, has the structure of the point-truncation of a Neumaier (sporadic $A_7$) geometry (similarly for $B \in \mathcal{B}$).
\end{enumerate}

\begin{rem} The Neumaier geometry is defined on a set $X$ of seven points, and has $35$ lines and $15$ planes. The set of lines consists of all $3$ element subsets of $X$. There exist $30$ different collections of $7$ lines yielding the structure of a projective plane. Under the action of the alternating group $A_7$, these projective planes fall into two orbits, each of size $15$. Choose one of these orbits as the planes of the Neumaier geometry; see \cite{cam}.
\end{rem}

\section{A geometric description of the Bouc complex of $McL$}

Recall that a $p$-subgroup $R$ of $G$ is called $p$-radical if $R = O_p(N_G(R))$. Let $\mathcal{B}_p(G)$ denote the collection of all non-trivial $p$-radical subgroups of $G$; this collection is known in the literature as the Bouc collection.

\medskip
In this section $G=McL$ and $\Delta$ represents the $2$-local geometry described in the previous section; further let $\mathcal{B} _2 = \mathcal{B} _2(McL)$. The Bouc complex of $McL$ was first determined in \cite{m} and later given by Yoshiara\footnote{This is from Table 5 in \cite{y05}. Notice that the table has a typographical error; for the third group $R=E_1^{(2)}$, $N(R)/R$ should read $(3 \times 3):2$ and not $S_3 \times S_3$.} in \cite{y05}. There are eight conjugacy classes of $2$-radical subgroups, as listed in {\it Table I}.

\medskip
\begin{center}
\begin{tabular}{l l l l l l}
{\bf Simplex} \phantom{a}& {\bf Name} \phantom{aa}& $\mathbf{R}$ \phantom{aaaa}& $\mathbf{Z(R)}$ \phantom{aa}&
$\mathbf{N_G(R)}$\phantom{aaaa}&{\bf Generators}\\
$p$&$R_p$ & $2$ & $A^1$ & $ 2^.A_8$&$\langle ab \rangle$\\
$A$&$R_A$ & $2^4$ & $A^{15}$ & $2^4:A_7$&$\langle a, b, c, d\rangle$\\
$B$&$R_B$ & $2^4$ & $A^{15}$ & $2^4:A_7$&$\langle a, b, e, f\rangle$\\
$AB$&$R_{AB}$ & $2^2.2^4$ & $A^3$ & $2^{2+4}(3 \times 3):2$&$\langle a, b, c, d, e, f\rangle$\\
$p\square$&$R_{p \square}$ & $2_+^{1+4}$ & $A^1$ & $2_+^{1+4}(S_3 \times S_3)$&$\langle a, b, c, e, u\rangle$\\
$pA\square$&$R_{pA\square}$ & $2.[2^5]$ & $A^1$ & $2.[2^5]:S_3$&$\langle a, b, c, d, e, u\rangle$\\
$pB\square$&$R_{pB\square}$ & $2.[2^5]$ & $A^1$ & $2.[2^5]:S_3$&$\langle a, b, c, e, f, u\rangle$\\
$pAB\square$&$R_{pAB\square}$ & $2^4:D_8$ & $A^1$ & $2^4:D_8$&$\langle a, b, c, d, e, f, u\rangle$\\
\end{tabular}

\vspace*{.5cm}
{\it Table I}
\end{center}

\medskip
The notation $\square$ stands for a  typical vertex of the truncated type in $\text{Res}_{\Delta}v$, where $v \in \lbrace \mathcal{P}, \mathcal{A}, \mathcal{B} \rbrace$. It corresponds to the square node
in the diagram of $\Delta$. Note that $\square$ does not represent a class of objects in the geometry $\Delta$. In $\text{Res}_{\Delta}p, \; \square$ stands for a line in the $PG(3,2)$ whose ``points" and ``planes" are the $\mathcal{A}$-spaces and $\mathcal{B}$-spaces incident with $p$. In $\text{Res}_{\Delta}A, \; \square$ can be identified with a point in the Neumaier geometry. The same is true when considering $\text{Res}_{\Delta}B$.

\medskip
In the last column of the {\it Table I} we give generators for each of the $2$-radical subgroups. Following Stroth \cite{str}, a presentation for the $2$-Sylow subgroup of $G$ can be written as $S = \langle a, b, c, d, e, f, u \mid \mathcal{R} \rangle$ with a set of relations:
\begin{align*}
\mathcal{R}  =& \lbrace a^2 = b^2 = c^2 = d^2 = e^2 = f^2 = u^2 =1, [a,b] = [a,c] = [a,d]=1\\
&[b,c]= [b,d] = [c,d]= [a, e]= [b,e]= [a,f]= [b,f]=[e,f]=1,\\
&c^e=abc, d^e=bd, c^f=bc, d^f=ad, a^u=b, c^u=c, d^u=cd, e^u=e, f^u=ef \rbrace\\
\end{align*}

In what follows we shall describe each of the subgroups of $\mathcal{B}_2$ in geometric terms.

\medskip
\hspace*{.5cm}{$R_{pAB\square} \simeq 2^4:D_8$}\\
The group $R_{pAB\square}$ of order $2^7$ is the $2$-Sylow subgroup of $G$. Its structure was analyzed in detail by Adem and Milgram in \cite[Lemma 1.1]{am97}.

\begin{lem}a) The group $R_{pAB\square}$ contains exactly one copy (a conjugate in $McL$) of each of the subgroups $R_A, R_B, R_{AB}$.\\
b)  The group $R_{pAB\square}$ contains exactly one copy of $R_{p\square}$.
\end{lem}

\begin{proof}a) See \cite{j68} for a proof.\\
b) We consider the composition of maps $2_+^{1+4} \rightarrow 2^4:D_8 \twoheadrightarrow UT_4 (2)$, where $UT_4 (2)$ denotes the group of upper triangular $4 \times 4$ matrices over ${\mathbf F}_2$. If this map is not injective, the central involution of $2_+^{1+4}$ lies in the center of $R_{pAB\square}$. The quotient $UT_4 (2)$ contains a unique subgroup of the form $2^4$ and its inverse image in $R_{pAB\square}$ is indeed a copy of $2^{1+4}_+$. If the above composition is injective, note that $UT_4 (2)$ contains a unique copy of $2_+^{1+4}$ and the inverse image of this subgroup in $R_{pAB\square}$ is of the form $2^{2+4}$. We reached a contradiction since $2^{2+4}$ does not contain a subgroup of the form $2_+^{1+4}$.
\end{proof}

There are $35$ involutions in $R_{pAB\square}$ and they can be partitioned as: $1+2+4_A+4_B+8_A+8_B+8_{\square}$, where $1$ stands for the central involution, $1+2$ represent the involutions contained in the second term of the upper central series of $R_{pAB\square}$ and $1+2+4_A+4_B$ the involutions contained in the third term of the same series. The involutions $1+2+4_A+8_A$ lie in the unique $R_A$ contained in $R_{pAB\square}$. Similarily $1+2+4_B+8_B$ lie in $R_B$. Note that $1+2$ correspond to the subgroup $2^2\simeq R_A \cap R_B$. Then $8_{\square}$ stands for the remaining $35-27$ involutions in $R_{pAB\square}$.

\medskip
\hspace*{.5cm}{$R_p \simeq 2$}\\
This corresponds to a single a central involution and thus a point of the geometry.

\medskip
\hspace*{.5cm}{$R_A \simeq 2^4$} \; \; and \; \;{$R_B \simeq 2^4$}\\
This is an elementary abelian $2$-subgroup of rank $4$. Its $15$ involutions are in one-to-one correspondence with the planes of the Neumaier geometry.

\medskip
\hspace*{.5cm}{$R_ {AB}\simeq 2^{2+4}$}\\
There are exactly $27$ involutions in this group. This structure corresponds to a pair $\lbrace A,B \rbrace$ of incident spaces from different classes. Furthermore $R_{AB} = \langle R_A ,R_B \rangle$, see \cite{fin73}. Clearly $R_A \subseteq R_{AB}$ and $R_B \subseteq R_{AB}$ are flags in $\mathcal{B}_2$. The involutions of $R_{AB}$ can be written as: $ 1+ 2 + 4_A + 8_A + 4_B + 8_B$.

\medskip
\hspace*{.5cm}{$R_ {p\square}\simeq 2_+^{1+4}$}\\
This group contains $19$ involutions which can be partitioned as: $1+2+4_A+4_B+8_{\square}$. Notice that $R_{p\square}$ and $R_{AB}$ are not incident in the Bouc complex. It is easy to see from the presentation given in {\it Table I} that $R_{p\square}$ is incident with both $R_{pA\square}$ and $R_{pB\square}$.

\medskip
\hspace*{.5cm}{$R_{pA\square} \simeq 2.[2^5]$} \; \; and \; \;{$R_{pB\square} \simeq 2.[2^5]$}\\
These groups contain $27$ involutions each. It is easy to see, from the presentations given in {\it Table I} that $R_A \subseteq R_{pA\square}$ and $R_B \subseteq R_{pB\square}$ are chains in $\mathcal{B}_2$. From this observation we deduce that the involutions of $R_{pA\square}$ can be written as: $ 1+ 2 + 4_A + 8_A + 4_B + 8_{\square}$ and similarily for $R_{pB\square}$: $ 1+ 2 + 4_A + 4_B + 8_B + 8_{\square}$.

\section{The homotopy relation}

In this Section we will prove:

\begin{thm} The Bouc complex $\mathcal{B}_2$ is homotopy equivalent to the subcomplex $\mathcal{B}_2^I$ obtained by removing from $\mathcal{B}_2$ the vertices of type $R_{p\square}, R_{pA\square}, R_{pB\square}$ and $R_{pAB \square}$ together with all the simplices containing them.
\end{thm}

{\it{Flags in $\mathcal{B}_2$}}. Let $F: \; x_1 \subseteq x_2 \ldots \subseteq x_n$ be a flag in $\mathcal{B}_2$. Then the stabilizer of the object $x_n$ acts on the collection of objects of type $x_1$ incident with $x_n$. In most of the cases this action is not transitive. The flags of the Bouc complex are given in {\it Table II}. In the first column of the table we give the type of the flag $F$, in the second column the order of the stabilizer of the object $x_n$ in $F$ and in the third column we give the number of elements in the orbits of $\text{Stab}_G(x_n)$ on the collection of objects of type $1$ incident with $x_n$.

\begin{proof}[{\bf Proof of the Theorem 4.1}]The proof of Theorem is straightforward, involving a sequence of steps of homotopy retractions which use the collapsibility property of certain flags. We will make repeated use of the following:

\begin{lem}[\cite{rsy}] Let $\Sigma \in \Delta$ be a simplex
of maximal dimension with $\sigma$ as a face. Assume that $\Sigma$ is the
only simplex of maximal dimension with $\sigma$ as a face. Then the process
of removing $\Sigma$ from $\Delta$, by collapsing $\Sigma$ down onto its
faces other than $\sigma$, is a homotopy equivalence.
\end{lem}

\medskip
\begin{scriptsize}
\begin{center}
\begin{tabular}{|l|l|l||l|l|l|}

\hline
&&&&&\\
{\bf Rank 1 flags}&&&{\bf Rank 3 flags}&& \\

$p$ & $2^7 \cdot 3^2 \cdot 5 \cdot 7$ & $1$ & $p \subseteq A \subseteq pAB \square$ & $2^7$ & $1,2,4,8$\\

$A$ & $2^7 \cdot 3^2 \cdot 5 \cdot 7$ & $1$ & $p \subseteq B \subseteq p AB \square$ & $2^7$ & $1,2,4,8$\\

$B$ & $2^7 \cdot 3^2 \cdot 5 \cdot 7$ & $1$ & $p \subseteq p \square \subseteq pAB \square$ & $2^7$ & $1,2,4,4,8$ \\

$pA\square$ & $2^7 \cdot 3$ & $1$ &  $p \subseteq pA\square \subseteq pAB \square$ & $2^7$
& $1,2,4,4,8,8$ \\

$pB\square$&$2^7 \cdot 3$ &$1$& $p \subseteq pB\square \subseteq pAB \square$ & $2^7$ & $1,2,4,4,8,8$\\

$AB$&$2^7 \cdot 3^2$ &$1$& $p \subseteq AB  \subseteq p AB \square$ & $2^7$ & $1,2,4,4,8,8$ \\

$p\square$ & $2^7 \cdot 3^2$ & $1$ &$p \subseteq A  \subseteq AB$ & $2^7 \cdot 3^2$ & $3,12$ \\

$p AB \square$& $2^7$&$1$&$p \subseteq B \subseteq AB$ & $2^7 \cdot 3^2$ &
$3,12$\\

{\bf Rank 2 flags}&&&$p \subseteq p\square \subseteq pA \square$ & $2^7 \cdot 3$ &
$1,6,12$ \\

$p \subseteq A$ & $2^7\cdot 3^2\cdot 5\cdot 7$ & $15$ &$p \subseteq p \square \subseteq pB \square$ &
$2^7 \cdot 3$ & $1,6,12$ \\

$p \subseteq B$ & $2^7 \cdot 3^2 \cdot 5 \cdot 7$ & $15$ &$p \subseteq A \subseteq p A \square$ & $2^7 \cdot 3$ & $1,6,8$ \\

$p \subseteq p \square$ & $2^7 \cdot 3^2$ & $1,18$ &$p \subseteq B \subseteq pB\square$ & $2^7 \cdot 3$ & $1,6,8$ \\

$p \subseteq AB$& $2^7 \cdot 3^2$& $3, 12, 12$&$A \subseteq AB \subseteq pAB\square$ & $2^7$ &
$1$ \\

$p  \subseteq pA\square$ & $2^7 \cdot 3$ & $1,6,8,12$
&$B \subseteq AB \subseteq pAB\square$ &
$2^7$ & $1$ \\

$p \subseteq pB\square$ & $2^7\cdot 3$ & $1,6,8,12$ &$A\subseteq pA \square \subseteq pAB\square$ & $2^7$ & $1$ \\

$p \subseteq pAB \square$ & $2^7$ & $1,2,4,4,8,8,8$  &$B\subseteq p B\square \subseteq pAB\square$
& $2^7$ & $1$ \\

$A  \subseteq pA\square$ & $2^7 \cdot 3$ & $1$ &$p \square \subseteq pA \square \subseteq pAB\square$ &
$2^7$ & $1$ \\

$A \subseteq AB$ & $2^7 \cdot 3^2$ & $1$ & $p\square \subseteq pB \square \subseteq pAB\square$&$2^7$&$1$\\

$A \subseteq pAB\square$ & $2^7$ & $1$ &&&\\

$B \subseteq pB\square$ & $2^7\cdot 3$ & $1$  &{\bf Rank 4 flags} && \\

$B \subseteq AB$&$2^7 \cdot 3^2$&$1$&$p \subseteq A \subseteq AB \subseteq pAB \square$   &$2^7$&$1,2,4,8$\\

$B \subseteq pAB\square$ &$2^7$&$1$&$p \subseteq B \subseteq AB \subseteq pAB \square$ &
$2^7$ & $1,2,4,8$\\

$p\square \subseteq pA\square$&$2^7\cdot 3$&$1$&$p \subseteq A \subseteq pA\square \subseteq pAB \square$ &
$2^7$ & $1,2,4,8$\\

$p\square \subseteq pB\square$&$2^7\cdot 3$&$1$&$p \subseteq B \subseteq pB\square \subseteq pAB \square$ & $2^7$ & $1,2,4,8$\\

$p\square \subseteq pAB\square$&$2^7$&$1$&$p \subseteq p\square \subseteq pA\square \subseteq pAB \square$ & $2^7$ & $1,2,4,4,8$\\

$AB \subseteq pAB\square$&$2^7$&$1$& $p \subseteq p\square \subseteq pB\square \subseteq pAB \square$& $2^7$ & $1,2,4,4,8$\\

$pA\square \subseteq pAB\square$&$2^7$&$1$&&&\\

$pB\square \subseteq pAB\square$&$2^7$&$1$&& &\\
&&&& &\\

\hline
\end{tabular}
\end{center}
\end{scriptsize}

\vspace*{.5cm}
\begin{center}
{\it Table II}
\end{center}

In what follows we will homotopically retract $\mathcal{B}_2$ down to $\mathcal{B}_2^I$. Recall that a maximal simplex $\Sigma$ is {\it free} over some maximal face $\sigma$, if $\Sigma$ is the only maximal simplex with $\sigma$ as a face. In this case we can use Lemma $3.1$ to collapse $\Sigma$ down onto its faces other than $\sigma$. This means that we can remove $\Sigma$ and $\sigma$ from the complex to obtain a homotopy equivalent subcomplex. Note that since this method depends on the maximality of $\Sigma$ and since the removal of some simplices might make other simplices maximal, the sequence of homotopy retractions we are performing has to be done in the order indicated here and listed in {\it Table III}.

\medskip
In the first column of the {\it Table III} we label the step. In the second column we give the maximal simplex $\Sigma$ and in the third column we give a face $\sigma$ of $\Sigma$ such that $\Sigma$ is free over $\sigma$. In the fourth column we specify the collection of points in $\sigma$.

\medskip
We give details for Step $6.A$ which is suggestive for the approach we use. Consider the simplex $\sigma : p \subseteq A \subseteq pA\square$. We are concerned here with the points denoted $2+4_A$. There are $6$ such points, $2$ points in each of the three flags of the form $p \subseteq A \subseteq pA\square \subseteq pAB\square$ containing $\sigma$. The element of order $3$ in $N_G(R_{pA\square})$ permutes the three Sylow subgroups containing $R_{pA\square}$, stabilizing the collection of involutions of $R_{pA\square}$; in particular, it permutes the three lines on the point $p$ in the plane $1+2+4_A$, thus we can write: $2+4_A = 2+2+2$. This means that for a fixed $\sigma$, there are $3$ simplexes $\Sigma$ of the type chosen above. We write this as $3 \cdot 2 \; p \subseteq  A \subseteq pA \square \subseteq pAB\square$, emphasizing also the $2$ points of $\Sigma$. Fix one such simplex $\Sigma$. This is a simplex of maximal dimension and free over the face $\sigma$, but only after the homotopy retractions of steps $1 - 5$. Thus we can apply Lemma $3.2$ and remove these two simplices. A similar argument can be used for the steps $6.B, 7, 12, 13$.

\medskip
After the 15 steps from {\it Table III} we are left with the complex $\mathcal{B}_2^I$, with a cone attached at each point $p$. This is of course homotopy equivalent to $\mathcal{B}_2^I$.
\end{proof}

\begin{rem}Note that the subcomplex $\mathcal{B}_2^I$ is essentially $| \mathcal{A}_2^- |$ from \cite[Section 6.11]{bs}. Further $\Delta$ corresponds to the semi-order complex $\Delta_{\mathcal{H}_I}$.  Note, also that $| \mathcal{A}_2^- |$ is a subdivision of $\Delta_{\mathcal{H}_I}$ where we subdivided by adding a vertex $AB$ (denoted $44'$ in \cite{bs}) and a new edge $p \subseteq AB$. It was proved in \cite{sy} that $\Delta$, $|\mathcal{A}_2^-|$ and the Quillen complex are homotopy equivalent.
\end{rem}

\begin{rem}
Using the information from {\it Table II}, the reduced Euler characteristic of the Bouc complex can be easily computed\footnote{This is given also in Table I of \cite{sy}.} and it is:
$$\tilde{ \chi} (\mathcal{B}_2) = 2^7 \cdot 7483$$
A $p$-subgroup $P$ of $G$ is called $p$-centric if $Z(P)$ is a Sylow $p$-subgroup of $C_G(P)$. Then $\mathcal{B}_p^{\rm{cen}}$ will denote the collection of $p$-radical and $p$-centric subgroups of $G$. Among the $2$-radical subgroups of McL all but the first conjugacy class of groups are $2$-centric.
The subcomplex $\mathcal{B}_2^{cen}$ of the centric subgroups of the group McL is the complex obtained by removing from $\mathcal{B}_2$ the vertex $R_p$ and all the flags containing it. This subcomplex has reduced Euler characteristic:
$$\tilde{ \chi} (\mathcal{B}_2^{cen}) = 2^6 \cdot 37241$$
Therefore, in this case $\Delta$ and $\mathcal{B}_2^{cen}$ are not homotopy equivalent.
\end{rem}

\bigskip
\begin{center}
\begin{tabular}{|l|l|l|l|}
\hline
Step&$\Sigma$&$\sigma$&Points in $\sigma$\\
\hline

$1.A$ & $p \subseteq A \subseteq AB  \subseteq pAB \square$ & $p \subseteq AB \subseteq pAB \square$ & $4_A+8_A$\\

$1.B$ & $p \subseteq B \subseteq AB  \subseteq pAB \square$ & $p \subseteq AB \subseteq pAB \square$ & $4_B+8_B$\\

$2.A$ & $p \subseteq A \subseteq pA\square \subseteq pAB\square$ & $p \subseteq pA\square \subseteq pAB\square$ &$8_A$\\

$2.B$ & $p \subseteq B \subseteq pB\square \subseteq pAB\square$ & $p \subseteq pB\square \subseteq pAB\square$ &$8_B$\\

$3.A$ & $p \subseteq p\square \subseteq pA\square \subseteq pAB\square$ & $p \subseteq pA\square \subseteq pAB\square$ & $4_B+8_{\square}$\\

$3.B$ & $p \subseteq p\square \subseteq pB \square \subseteq pAB\square$ & $p \subseteq  pB\square \subseteq pAB\square$ &$4_A+8_{\square}$\\

$4.A$ & $p \subseteq  A \subseteq pA \square \subseteq pAB\square$ & $p \subseteq A \subseteq pAB\square$ &$4_A$\\

$4.B$ & $p \subseteq  B \subseteq pB \square \subseteq pAB\square$ & $p \subseteq B \subseteq pAB\square$ &$4_B$\\

$5.A$ & $p \subseteq p\square \subseteq  pA\square \subseteq pAB\square$ & $p \subseteq pA\square \subseteq pAB\square$ &$4_A$\\

$5.B$ & $p \subseteq p\square \subseteq  pB\square \subseteq pAB\square$ & $p \subseteq pB\square \subseteq pAB \square$ &$4_B$\\

$6.A$ & $3\cdot 2p \subseteq A \subseteq pA\square \subseteq pAB\square$ & $p \subseteq A \subseteq pA\square$ &$2+4_A$\\

$6.B$ & $3\cdot 2p \subseteq B \subseteq pB\square \subseteq pAB\square$ & $p \subseteq B \subseteq pB\square$ &$2+4_B$\\

$7.A$ & $3\cdot 2p \subseteq p\square \subseteq pA\square \subseteq pAB\square$ & $p \subseteq p\square \subseteq pA\square$ &$2+4_A$\\

$7.B$ & $3\cdot 2p \subseteq p\square \subseteq pB\square \subseteq pAB\square$ & $p \subseteq p\square \subseteq pB\square$ &$2+4_B$\\

$8.A$& $p \subseteq A \subseteq AB \subseteq pAB\square$ & $p \subseteq A \subseteq pAB\square$ &$2$\\

$8.B$& $p \subseteq B \subseteq AB \subseteq pAB\square$ & $p \subseteq B \subseteq pAB\square$ &$2$\\

$9.A$ & $p \subseteq A \subseteq pAB\square$ & $p \subseteq pAB\square$ & $8_A$\\

$9.B$ & $p \subseteq B \subseteq pAB\square$ & $p \subseteq pAB\square$ & $8_B$\\

$10$ & $p \subseteq p\square \subseteq pAB\square$ & $p \subseteq pAB\square$ & $4_A+4_B+8_{\square}$\\

$11.A$ & $p \subseteq p\square \subseteq pA\square$ & $p \subseteq pA\square$&$4_B+8_{\square}$\\

$11.B$ & $p \subseteq p\square \subseteq pB\square$ & $p \subseteq pB\square$&$4_A+8_{\square}$\\

$12$ & $9 \cdot 2 p \subseteq p\square \subseteq pAB\square$ & $p \subseteq p\square$ & $2+4_A+4_B+8_{\square}$\\

$13.A$ & $3\cdot 2 p \subseteq pA\square \subseteq pAB\square$ & $ p \subseteq pA\square$&$2+4_A$ \\

$13.B$ & $3\cdot 2 p \subseteq pB\square \subseteq pAB\square$ & $ p \subseteq pB\square$&$2+4_B$ \\

$14.A$ & $p \subseteq A \subseteq pA\square$ & $p \subseteq pA\square$ &$8_A$\\

$14.B$ & $p \subseteq B \subseteq pB\square$ & $p \subseteq pB\square$ &$8_B$\\

$15$ & $p \subseteq AB \subseteq pAB\square$ & $p \subseteq pAB\square$ &$2$\\

\hline
\end{tabular}

\vspace*{.5cm}
{\it Table III}
\end{center}


\medskip
\begin{thebibliography}{20}

\bibitem{am97}A. Adem, R.J. Milgram, {\it The cohomology of the McLaughlin group and some associated groups}, Math. Z. {\bf 224} (1997), 495-517.

\bibitem{bs}D. Benson, S.D. Smith, {\it Classifying spaces of sporadic groups}, preprint available at http://www.maths.abdn.ac.uk/$\sim$bensondj/html/archive/benson-smith.html, 2006.

\bibitem{cam}P.J. Cameron, {\it Projective and polar spaces}, QMW Math. Notes {\bf 13}, 1991.

\bibitem{fin73}L. Finkelstein, {\it The maximal subgroups of Conway's group $C_3$ and McLaughlin's group}, J. Algebra {\bf 25} (1973), 58-89.

\bibitem{j68}Z. Janko, {\it Characterization of the Mathieu simple groups}, J.
Algebra {\bf 9} (1968), 20-41.

\bibitem{mgo1}J.S. Maginnis, S.E. Onofrei, {\it On a homotopy
relation between the
$2$-local geometry and the Bouc complex for the sporadic group Co$_3$}, arXiv: math.GR/0510542, 2005.

\bibitem{mgo2}J.S. Maginnis, S.E. Onofrei, {\it New collections of $p$-subgroups and homology decompositions for classifying spaces of finite groups}, arXiv: math.GR/0510540, 2005.

\bibitem{m}J.C. Murray, {\it Dade's conjecture for the McLaughlin simple group}, PhD Thesis, University of Illinois at Urbana-Champaign, 1998.

\bibitem{rst}M.A. Ronan, G. Stroth, {\it Minimal parabolic geometries for the
sporadic groups}, Europ. J. Combin. {\bf 5} (1984), 59-91.

\bibitem{rsy}A.J.E. Ryba, S.D. Smith, S. Yoshiara, {\it Some projective modules
determined by sporadic geometries}, J. Algebra {\bf 129} (1990), 279-311.

\bibitem{sy}S.D. Smith, S. Yoshiara, {\it Some homotopy equivalences for sporadic
geometries}, J. Algebra {\bf 192} (1997), 326-379.

\bibitem{str}G. Stroth, {\it Some geometry for McL}, Comm. Algebra {\bf 17}(11) (1989), 2825-2833.

\bibitem{y05}S. Yoshiara, {\it Radical 2-subgroups of the Monster and Baby Monster}, J. Algebra {\bf 287} (2005), 123-139.\\
\end{thebibliography}
\end{document}